\documentclass[10pt]{article}

\usepackage{amssymb,amsmath}
\usepackage{url}
\usepackage{graphicx}
\usepackage{tikz}
\usepackage{algpseudocode}
\usepackage{multirow}
\usepackage{placeins}
\usepackage{subfig}
\usepackage{enumerate}

\newcommand{\tbbnosp}{\textsc{TBB}}
\newcommand{\tbb}{\tbbnosp\ }
\newcommand{\OpenMPReg}{OpenMP}
\newcommand{\inteltbb}{Intel Threading Building Blocks}
\newcommand{\ds}{\displaystyle}
\def\tlvs{\vrule height 1em  width 0pt} 

\title{Task-based adaptive multiresolution for time-space multi-scale
  reaction-diffusion systems on multi-core architectures}
  
\author{
St\'ephane~Descombes\footnotemark[1]
\and 
Max~Duarte\footnotemark[2]
\and 
Thierry~Dumont\footnotemark[3]
\and 
Thomas~Guillet\footnotemark[4]
\and 
Violaine~Louvet\footnotemark[3]
\and 
Marc~Massot\footnotemark[5]
}

\begin{document}

\maketitle

\renewcommand{\thefootnote}{\fnsymbol{footnote}}

\footnotetext[1]{
Universit\'e C\^ote d'Azur, CNRS, Inria, LJAD,  France
({\tt stephane.descombes@unice.fr}).
}
\footnotetext[2]{
CCSE,
Lawrence Berkeley National Laboratory,
1 Cyclotron Rd. MS 50A-1148,
Berkeley, CA 94720, USA
and
CD-adapco, 200 Shepherds Bush Road, London W6 7NL, UK
({\tt max.duarte@cd-adapco.com}).
}
\footnotetext[3]{
Universit\'e de Lyon, Universit\'e Claude Bernard Lyon 1,
CNRS UMR 5208, Institut Camille Jordan, 43 blvd. du 11 novembre 1918, 
F-69622 Villeurbanne Cedex,
France 
({\tt tdumont@math.univ-lyon1.fr}, {\tt louvet@math.univ-lyon1.fr}).
}
\footnotetext[4]{
Intel,
Les Montalets,
2 rue de Paris,
92196 Meudon,
France
and
Exascale Computing Research,
Campus Teratec,
2 rue de la Piquetterie,
91680 Bruy\`eres-le-Ch\^atel,
France
({\tt thomas.guillet@intel.com}).
}
\footnotetext[5]{
CNRS UPR 288, Laboratoire EM2C, 
CentraleSup\'elec,
F\'ed\'eration de Math\'ematiques de l'\'Ecole Centrale Paris, CNRS FR 3487, Grande Voie des Vignes, 92295 Chatenay-Malabry Cedex,
France
({\tt marc.massot@centralesupelec.fr}).
}

\renewcommand{\thefootnote}{\arabic{footnote}}

\begin{abstract}
A new solver featuring time-space adaptation and error control has been recently introduced
to tackle the numerical solution of stiff reaction-diffusion systems. 
Based on operator splitting, finite volume adaptive multiresolution 
and high order time integrators with specific stability properties for each operator, 
this strategy yields high computational efficiency 
for large multidimensional computations on standard architectures such as powerful workstations.
However, the data structure of the original implementation, 
based on trees of pointers, provides limited opportunities for efficiency enhancements,
while posing serious challenges in terms of parallel programming and load balancing.
The present contribution proposes a new implementation of the whole set of numerical methods including \textsc{Radau5} and \textsc{ROCK4},
relying on a fully different data structure together with the use of a specific library, \textsc{TBB}, for 
shared-memory, task-based parallelism with work-stealing.
The performance of our implementation is assessed in a series of test-cases of increasing difficulty in 
two and three dimensions on multi-core and many-core architectures, demonstrating high scalability.
\end{abstract}

\noindent {\bf Keywords:}
Task-based parallelism, multi-core architectures, multiresolution, adaptive grid, stiff reaction-diffusion equations

\noindent {\bf AMS subject classifications:}
65Y05, 65T60, 65M50, 65L04, 35K57

\pagestyle{myheadings}
\thispagestyle{plain}
\markboth{DESCOMBES, DUARTE, DUMONT, GUILLET, LOUVET, MASSOT}
{TASK-BASED ADAPTIVE MULTIRESOLUTION ON MULTI-CORE ARCHITECTURES}

\section{Introduction}
Stiff reaction-diffusion systems 
model various complex phenomena across different disciplines such as combustion, atmospheric sciences, plasma physics or biomedical engineering. 
Such models can also involve the dynamics of moving fronts, usually very localized in space (see, \emph{e.g.}, \cite{duarte:tel-00667857} and references therein), 
and potentially entail a large number of unknowns. 
A general reaction-diffusion system can be written as follows,
for $i=1,2,\ldots,m$:
\begin{equation}\label{thesystem}
\left.
\begin{array}{ll}
\displaystyle
\frac {\partial u_i}{\partial t}(x,t)- {\rm{div}}( \varepsilon_i(x)
\ {\rm{grad }}\ u_i (x,t))= f_i(u(x,t)),
&
x \in \Omega \subset \mathbb{R}^d,\, t>0,\\[2ex]
u_i(x,0)= u_i^0(x),
&
x \in \Omega;
\end{array}
\right\}
\end{equation}
with the compact notation: $u=(u_1,\ldots, u_m)^t$.
In particular we denote $f(u)=(f_1(u),\ldots, f_m(u))^t$.
With no loss of generality 
we restrict our presentation to homogeneous Neumann boundary
conditions.

Two major difficulties need to be addressed
when solving numerically this kind of problem in two and three dimensions.
First, 
a large spectrum of temporal scales in the nonlinear source terms yields 
highly stiff equations\footnote{The latter can be expected when the Jacobian
matrices, $(\partial f_i / \partial u_j)_{1 \leq i,j \leq m}$, have
eigenvalues whose real part varies within a large negative interval or in the presence of strong attracting slow manifolds \cite{MR1439506}.}.
Systems of  stiff ordinary differential equations
impose the use of numerical methods with specific stability properties in order to achieve accuracy and stability within reasonable memory and computing costs \cite{MR1439506}.
Secondly,
steep fronts require a very fine discretization mesh, at least locally, which leads to problems of large size if no mesh adaptation is used. 
Additionally, spatial stiffness may arise 
as a consequence of these steep spatial gradients
even with non-stiff source terms and 
diffusion coefficients of relatively small value \cite{MR2335366}.

We have recently introduced in \cite{MR2890259}
a tailored numerical strategy to cope with the latter difficulties using reasonable computing resources, that is, 
on a sufficiently powerful workstation, possibly exploiting shared-memory parallelism \cite{MR2869529}.
It relies on 
time operator splitting\footnote{High performance computing can be achieved 
by choosing dedicated schemes for each split sub-system
\cite{MR2890259,Dumont2013}. 
Operator splitting schemes also
exhibit a large data-driven parallelism,
and their mathematical analysis is 
well-established for relatively large splitting time steps
even when stiffness is present
\cite{MR2127928,MR2335366,MR3213800}.}
with one-step, high order integration schemes, namely,
\textsc{Radau5}\footnote{An implicit, fifth order Runge-Kutta scheme with $A$- and $L$-stability properties and time step adaptation.}
\cite{MR1439506} and
\textsc{Rock4}\footnote{A stabilized, explicit Runge-Kutta method of order four and time step adaptation.} 
\cite{MR1923724}
for the time integration of the reaction and diffusion sub-systems, 
respectively.
The discretized equations are solved on a dynamically adapted grid generated by 
multiresolution analysis (MRA) in a finite volume multi-dimensional framework in the spirit of the original 
work of Harten \cite{Harten95} and then  Cohen et al. and M\"{u}ller
\cite{Cohen03,muller2003,MR2647599}. MRA is based on a wavelet decomposition and a multiresolution transform, yielding both
 highly compressed representations for problems displaying localized fronts as well as a compression error control 
 with respect to the solution on the full grid. Let us emphasize that this property is a major difference with an AMR (Adaptive Mesh Refinement) strategy 
 where the refinement criterion rather relies on heuristics. 
This MRA-operator splitting numerical strategy  was implemented 
in the \textsc{MBARETE} code \cite{duarte:tel-00667857},
using a tree-structured data with pointers and recursive navigation. 
However, the latter implementation shows serious 
limitations in terms of parallel programming due to the lack of data locality and load balancing.
A new paradigm of parallelism together with a different and customized data structure is thus
needed to achieve more efficient MRA-operator splitting implementations.

A wide spectrum of techniques and runtime implementations are available for application developers to express parallelism.
For the problems we are interested in and using MRA, one single modern compute node 
provides enough memory and computing power to run simulations with a reasonable time-to-solution.
Therefore, we have so far focused on parallelism over shared-memory architectures as it
provides a number of advantages well-adapted to MRA applications.
Contrary to uniform Cartesian grids for which arrays can be generally accessed following regular patterns within long loops,
adaptive meshing in general, that is both multiresolution and AMR codes, relies on fine-grained and dynamic data structures.
The corresponding algorithms can thus have intricate data dependencies, especially for complex operations such as those associated with mesh adaptation.
Hence, exposing parallelism in these methods is best done using programming techniques that combine both coarse- and fine-grained parallel constructs,
keeping in mind that maximizing parallel coverage is crucial to limiting the impact of Amdahl's law on scalability.
As a consequence of the complex layouts of the resulting adapted meshes, 
an efficient parallel implementation requires non-trivial balancing of the computations across all available computing cores.
Moreover, the load balancing needs to be frequently updated as the mesh structure evolves during the simulation.
Actually, many existing mesh refinement algorithms and libraries address this issue in distributed-memory settings for both patch- and cell-based AMR \cite{crutchfield_welcome_1993,oliker_plum_1998,MR1729306,rendleman_parallelization_2000,macneice_paramesh:_2000,wissink_large_2001,teyssier_cosmological_2002,deiterdingAMROC,burstedde_p4est:_2011,MR2869525,keppens_parallel_2012,drui2016, essadki2016}. 
However, even if MRA leads to compression error control, 
the mesh adaptation relies on a multiresolution transform that involves recursive navigations throughout the data sets \cite{Brix11}.
The latter makes it much more difficult to parallelize than classic AMR approaches. Some advances have been conducted recently into that direction but in a different framework
  \cite{forster2016,SpaceX_gputechconf,Nejadmalayeri2015237}.

In shared-memory architectures, the common memory address space enables dynamic load balancing while avoiding data transfers, but it also requires careful synchronizations between threads.
Starting from the original implementations 
on a single computer \cite{crutchfield_welcome_1993,MR1257158},
several authors have proposed multithreaded implementations for adaptive gridding techniques \cite{balsara_highly_2001,dreher_racoon:_2005}
with some focusing on specialized grid structures for multi-core processors \cite{hutchison_blocking_2010,hutchison_cluster_2013}.
Task-based parallelism provides an attractive shared-memory solution to both granularity and dynamic load balancing requirements.
In a task-based approach the programmer introduces parallelism by specifying computations that can be carried out in parallel.
Expressing tasks can be done using different techniques, for example, relying on compiler directives or through library calls.
Scheduling of the tasks is determined at runtime based on the available computing resources, yielding dynamic load balancing using techniques such as work-stealing \cite{blumofe1996cilk,Blumofe1999}.
A key feature associated with tasks is that parallelism is not limited to flat parallel iteration constructs, but can be introduced recursively by having tasks create other tasks.
This makes the task concept particularly suited to codes with complex hierarchical operations and recursive navigation in the data set, such as MRA applications.

We have thus opted for a shared-memory parallel approach based on work-stealing, relying
on the Intel Threading Building Blocks (\tbbnosp) library \cite{Reinders:2007:ITB:1461409,tbbwebsite}.
While other runtime libraries supporting task parallelism are currently available, such as the 
 {\OpenMPReg} application programming interface (API), and could have been used for our purposes with equivalent functionality, 
 they would have required a much more important programming effort.
Shared-memory parallelism was hence introduced using the task-based {\inteltbb} runtime library.
The latter turns out to be particularly well suited to the multiresolution algorithm in order to cope with nested
parallel regions, which are directly related to the MRA-splitting strategy.

Additionally, a completely new implementation of the original numerical strategy of \textsc{MBARETE} code has been conducted,
where a different and more efficient data structure has been introduced as an alternative to the original pointer-based approach.
Important performance enhancements due to work-stealing techniques and in particular the use of \tbb
have already been shown in \cite{Hejazialhosseini2010,Rossinelli11}
for multiresolution schemes applied to non-stiff time-dependent PDEs, that is, 
using explicit time integration schemes\footnote{These authors developed a
data structure based on the definition of {\it wavelet blocks},
specifically a tree of wavelet blocks, where each block contains a predefined number of cells at the same grid-level.}.

In order to assess the new implementation and code performance we have chosen three reaction-diffusion models
with increasing complexity in two ($d=2$) and three ($d=3$) dimensions:

\begin{enumerate}[I.]
\item {\bf NAGUMO.} A bistable, Nagumo-type reaction-diffusion equation\footnote{A bistable case of type A according 
to \cite{volpert94_3}, which is a limit case of the Nagumo function,
similar to what is found in combustion models \cite{zeldovich85} or in nonlinear chemical dynamics 
when studying traveling waves due, for instance, to chemical reactions with cubic auto-catalysis \cite{grayscott94}.}.
Here we have $m=1$ and $f_1(u_1)= k u_1^2 (1-u_1)$ in (\ref{thesystem}). 
We consider $k=10$, $\varepsilon_1= 0.1$, and an initial condition satisfying 
$0 \leq u^0_1(x) \leq 1$ everywhere\footnote{In general the reaction term is not stiff as 
$|df_1(u_1)/du_1| \leq 3k$.  However, with this setting, traveling waves develop with a sharp spatial gradient 
yielding a space multi-scale configuration.}.

\item {\bf BZ.} The Belusov-Zhabotinsky reaction \cite{noyes1972oscillations,jahnke1989chemical}. 
This is a system of three ($m=3$) equations, 
where the reaction term is given by\footnote{This set describes a
chemical reaction between $HBrO_2$, $Ce(IV)$ and $Br^{-}$, also known
as the Oregonator problem.}
$f_1(u_1,u_2,u_3)=10^5\ (-0.02\ u_1- u_2 u_1+ 1.6\  u_3)$,
$f_2(u_1,u_2,u_3)= 10^2\ (u_2 - u_2^2-  u_1(u_2-0.02))$, and
$f_3(u_1,u_2,u_3)=u_2 - u_3$.
For the diffusion part, we set $\varepsilon_1(x)=\varepsilon_2(x)=2.5\times 10^{-3}$
and $\varepsilon_3(x)=1.5\times 10^{-3}$ in (\ref{thesystem})\footnote{The reaction term is stiff.
The system of ordinary differential equations has a
limit cycle along which the amplitude of the eigenvalues of the Jacobian attains values of the order of $10^5$.
As a comparison, when computing in 
$\Omega= [0,1]^d$ with a spatial discretization step of $h=1/1024$,
the largest negative eigenvalues in the diffusion term are about $-2\times 10^4$.
Propagating fronts with steep spatial gradients are also developed in this case
yielding a time-space multi-scale configuration.}.

\item {\bf STROKE.} An ischemic stroke model \cite{stroke,ddstroke}. This is a system of 21 ($m=21$)
equations with a very stiff reaction term\footnote{The right-hand side, $f(u)$,
incorporates an important amount of biological knowledge about ion channels
(see \cite{Dumont2013} for the model here considered).
A numerical computation of the eigenvalues of the
Jacobian near a stable equilibrium state, $f(u)=0$, shows real parts in the
interval of $[-10^8,0[$. Stiffness is due in part to the fact that $f(u)$
models voltage-gated ion channels that open or close when the difference of potential
between a cell and the surrounding media attains a given threshold. 
Gates are modeled through sigmoid functions, closely approximating a
Heaviside function. Simulations show steep spatial gradients, as well.}.
Its computation is performed by a quite complex and computationally expensive
program of about $400$ instructions, 
containing many $\log$-function evaluations.
\end{enumerate}

The paper is organized as follows.
After a brief description of the key aspects of the numerical strategy in Section \ref{sec:strategy},
we focus in Section \ref{sec:implementation} on the task-based parallel implementation developed in this work.
A performance analysis is then presented in Section \ref{sec:performance},
to end with some concluding remarks and future developments related to this work.

\section{Numerical strategy}\label{sec:strategy}
We are interested in two- and  three-dimensional simulations. 
Discretizing (\ref{thesystem}) in space yields a large and stiff system of ordinary
differential equations given by
\begin{equation}\label{discretise}
 \frac{d U}{dt}=   A_{\varepsilon} U + F(U),
\end{equation}   
where $A_{\varepsilon}$ is a matrix corresponding to the discretization
of the diffusion operator. Using a standard centered second-order discretization 
on a Cartesian mesh, 
$A_{\varepsilon}$ is given by a matrix with five (resp., seven) non-zero elements
per line in two (resp., three) dimensions.
Following \cite{MR2890259} a dedicated splitting solver, briefly described in 
what follows, is implemented to 
advance system (\ref{discretise}) in time.

\subsection{Dedicated splitting solver}\label{sec:splitting}
Considering independently the two sub-problems coming from (\ref{discretise}):
\begin{equation*}
\begin{array}{lcr}
\left.
\begin{array}{l}
\displaystyle \frac{d V}{dt}= A_{\varepsilon}V,\\[1.75ex]
V(0)=V_0,
\end{array}
\right\}
&
\qquad
&
\left.
\begin{array}{c}
\displaystyle \frac{d W}{dt}= F(W),\\[1.75ex]
W(0)=W_0,
\end{array}
\right\}
\end{array}
\end{equation*}
we denote by $D_{\Delta t}V_0$ and $R_{\Delta t}W_0$
the solution of the first and second sub-problem,
respectively,
after a given \emph{splitting time step} $\Delta t$.
Assuming in general stiff reaction terms,
we consider without any loss of generality the 
second-order Strang scheme\footnote{It has been shown 
that better accuracies are expected by ending the
splitting scheme with the substep involving the fastest time scales \cite{MR2014215,MR2127928,MR3213800}.
Both theoretical and numerical results
show that this method performs very well for stiff reaction-diffusion
systems \cite{duarte:tel-00667857,MR2890259,Dumont2013}.}
\cite{MR0235754},
$U_{n+1}=R_{\Delta  t/2} \circ  D_{\Delta t}\circ R_{\Delta  t/2} U_{n}$.

The \textsc{Radau5} method is used to solve 
each of the independent initial value problems coming from the reaction term
and defined at each grid-point, 
during the \emph{Reaction Steps}, $R_{\Delta  t/2}$; whereas
the \emph{Diffusion Step}, $D_{\Delta  t}$, uses the explicit \textsc{ROCK4} method
to globally advance the discrete diffusion sub-problem.
Notice that with the latter choice of solver, one 
only needs to evaluate matrix-vector products without any need of
preconditioning nor factorizing the diffusion matrices when the mesh has changed.
In this work 
we have chosen splitting time steps, 
$\Delta t$, equal to $10^{-2}$, 
$10^{-3}$, and $1$ for
the NAGUMO, BZ, and STROKE models, respectively.

Finally, in terms of parallelism
the $R_{\Delta  t/2}$-step exhibits a high data-driven parallelism,
for one has as many independent problems as nodes in the spatial discretization.
The situation is different for the $D_{\Delta  t}$-step, which involves the whole
computational domain.
A simple and straightforward parallelism consists in parallelizing the solution of
the $m$ independent parabolic linear equations in (\ref{thesystem}) \cite{MR2869529}.
However, 
an alternative that is not limited by the number of equations
should be considered to enhance the parallel performance of the diffusion solver.
The latter is much simpler to achieve in the case of an explicit scheme like \textsc{ROCK4},
since matrix-vector products are easier to parallelize than 
linear solvers.

\subsection{Adaptive multiresolution method}
For the kind of problems we are interested in, featuring propagating localized fronts,
an adaptive meshing technique that reduces the number
of grid-nodes over near-equilibrium zones can significantly reduce
the CPU time\footnote{Since this is a key aspect, 
a numerical experimentation on uniform meshes is provided in Appendix \ref{fixedmeshes}
for the three aforementioned models.} together with the required computer memory.
In the following we describe some basic features of the multiresolution method 
introduced in \cite{Cohen03}, and implemented in \cite{MR2890259} 
with the aforementioned dedicated splitting solver
for reaction-diffusion problems\footnote{More details on this particular multiresolution implementation can be found in \cite{duarte:tel-00667857}.}.
\begin{figure}[h]
\begin{center}
\includegraphics[width=0.85\linewidth]{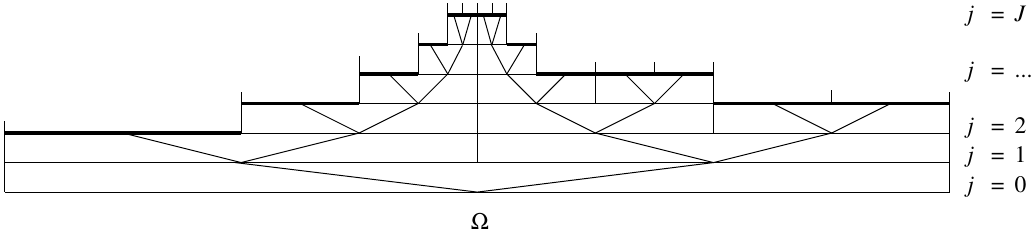}
\caption{One-dimensional graded tree composed of nested grids.
The resulting adapted mesh is given by the leaves of the tree, indicated in bold.}
\label{labelgrid}
\end{center}
\end{figure}

Without loss of generality we consider the computational domain 
$\Omega=[0,1]^d$, $d=1,2,3$. 
A recursive dyadic subdivision of $\Omega$ is then performed, 
yielding $2^{dj}$ cells (segments, squares or cubes) of equal size 
at each grid-level $j$, from $j=1$ to $j=J$. 
Levels $0$ and $J$ correspond, respectively, to the whole computational 
domain $\Omega$ and the mesh at the finest spatial resolution.
An adapted spatial representation can be thus achieved by 
trimming the set of nested grids throughout the tree, as shown in Figure~\ref{labelgrid} 
for a one-dimensional case.
In particular the adapted grid is given by the \emph{leaves} of the tree.
Finite volumes are associated with all cells of the tree, and throughout this paper
we will refer to all cells, including the leaves, and their corresponding finite volumes 
as \emph{nodes}.
The \emph{root} of the tree is the entire domain $\Omega$. 
Even though the solution of problem (\ref{thesystem}) is only performed
on the adapted mesh, that is on the leaf nodes, the solution is
updated and stored at all nodes of the tree. 
The set of values at nodes of level $j$ is denoted as $\mathbf{U}_j$.
Data at different levels of discretization are related
by two inter-level transformations defined by the \emph{projection}
and \emph{prediction} operators.
The projection operator, $P^j_{j-1}$, which maps $\mathbf{U}_j$ to
$\mathbf{U}_{j-1}$,
is obtained through exact averages computed at the finer level;
that is, for a given node at level $j-1$ one computes the 
average of the  $2^{d}$ nested nodes at the successive finer level $j$.
The prediction operator, $P^{j-1}_j$,
maps $\mathbf{U}_{j-1}$ to 
an approximation $\hat{\mathbf{U}}_j$ of $\mathbf{U}_{j}$;
that is, the value of a given node at level $j$ is approximated based on the surrounding coarser nodes 
at level $j-1$.

Grid adaptation is performed based on
local estimators of the spatial regularity of the solution at a given 
simulation time. 
These local estimators are known as \emph{details},
and for a given node at level $j$ its detail is defined as
\begin{equation}\label{eq:detail}
d_{j,k}=u_{j,k}-\hat{u}_{j,k} =
u_{j,k}-P^{j-1}_j\circ P^j_{j-1} u_{j,k},
\end{equation}
where $k \in \mathbb{Z}^d$ accounts for the location of the node at level $j$,
and $u_{j,k}$ represents the cell-average of $u(x,t)$ there.
From the theoretical point of view, the details correspond to the coefficients
related to a given discrete function when represented on a wavelet basis.
This wavelet decomposition is a general procedure, 
independent of any physical particularity of the problem.
Introducing a tolerance parameter, 
$\epsilon >0$, threshold values are defined level-wise as
$\epsilon_j =2^{\frac{d}{2}(j-J)}\epsilon$, $j=0,1,\ldots,J$.
Data compression is thus achieved by discarding nodes
whose details are smaller than $\epsilon_j$
in a given norm; 
here we consider an $L_2$-norm \cite{duarte:tel-00667857}.
Conversely, 
nodes whose details exceed
$\epsilon_j$ must be refined. 
Grid adaptation is hence driven by the local value of the details,
while the tolerance parameter defines the level of approximation errors
introduced by the data compression \cite{Cohen03}.
This is the main difference between MRA and AMR strategies. In particular, while enabling compression error control,
the multiresolution analysis also leads to recursive navigation in the data set.

In this work 
the approximated values, $\hat{\mathbf{U}}_j$, generated with the prediction operator,
are obtained using centered polynomial interpolations 
of accuracy order, $o=2l+1$, computed
with the $l$ nearest neighboring cells in each direction
including the diagonals in multi-dimensional configurations.
For a one-dimensional configuration with $l=1$,
the prediction operator is then given by \cite{Cohen03}
\begin{equation}
\hat{u}_{j+1,2k}=u_{j,k} + \frac{1}{8}(u_{j,k-1} - u_{j,k+1}),
\qquad 
\hat{u}_{j+1,2k+1}=u_{j,k} + \frac{1}{8}(u_{j,k+1} - u_{j,k-1}).
\label{inter_1d} 
\end{equation}
Extension to multi-dimensional Cartesian grids is easily
obtained by a tensor product of the one-dimensional operator (\ref{inter_1d})
(see, \emph{e.g.}, \cite{Roussel03}). 
The interpolation stencil is thus given by $(2l+1)^d$ cells at the coarser level.
Here we will only consider third-order interpolations with $l=1$.
Notice that 
these local interpolation stencils must be available throughout the tree-structure
during the simulation.
A tree-structure that complies with such a constraint is called a \emph{graded tree} \cite{Cohen03}.
Hence, a graded tree-structure 
must be always guaranteed after the refining and coarsening operations.

\section{Task-based parallelism for MRA-splitting applications}\label{sec:implementation}

We have introduced shared-memory parallelism using the task-based {\inteltbb} runtime library \cite{Reinders:2007:ITB:1461409,tbbwebsite}.  
It turns out to be particularly well suited to the multiresolution algorithm in order to cope with nested
parallel regions, which are directly related to the MRA-splitting strategy. 
Other than TBB, several runtime libraries supporting task parallelism are currently available, 
such as the {\OpenMPReg} API \cite{OpenMP}, XKaapi \cite{XKaapi}, among others.
They could potentially have been used for our purposes with equivalent functionality, but requiring a much more  important programming effort. 
A discussion clarifying the reasons why we have opted for \tbbnosp, in particular  over {\OpenMPReg} API, is provided in Appendix \ref{TBBvsOpenMP}.

\subsection{Parallel implementation of the multiresolution algorithm}

\begin{figure}[h]
\begin{center}
\scalebox{0.8}{\begin{tikzpicture}


\tikzstyle{stencil point}=[circle, inner sep=1mm, fill=black]

\begin{scope} 
	\node[left] at (0, 1) {$Q_1$};
	\draw[thick, dotted] (0, 0) grid +(2,2);
	\node[stencil point] at (0.5, 1.5) {};
	\node[stencil point] at (1.5, 1.5) {};
\end{scope}

\begin{scope}[yshift=2.2cm] 
	\node[left] at (0, 1) {$Q_2$};
	\draw[thick] (0, 0) grid +(2,2);
	\draw (1.5, 0) -- (1.5, 1);
	\draw (1, 0.5) -- (2, 0.5);
	\draw (1.25, 0.25) node {\small $P$};
	\draw (1.75, 0.25) node {\small $P$};
	\draw (1.75, 0.75) node {\small $P$};
	\draw (1.25, 0.75) node {\small $P$};
	\node[stencil point] at (0.5, 0.5) {};
	\node[stencil point] at (0.5, 1.5) {};
	\node[stencil point] at (1.5, 1.5) {};
	\node[stencil point] at (1.5, 0.5) {};
\end{scope}

\begin{scope}[xshift=2.2cm] 
	\node[right] at (2, 1) {$Q_3$};
	\draw[thick, dotted] (0, 0) grid +(2,2);
	\node[stencil point] at (0.5, 1.5) {};
\end{scope}

\begin{scope}[xshift=2.2cm, yshift=2.2cm] 
	\node[right] at (2, 1) {$Q_4$};
	\draw[thick, dotted] (0, 0) grid +(2,2);
	\node[stencil point] at (0.5, 0.5) {};
	\node[stencil point] at (0.5, 1.5) {};
\end{scope}

\draw (0.2, -0.7) node[stencil point] {} +(0.2, 0) node[right] {Stencil nodes};

\begin{scope}[
	xshift=10cm,
	grow'=up,
	level 1/.style={sibling distance=2.4cm},
	level 2/.style={sibling distance=0.5cm},
	level 3/.style={sibling distance=0.5cm, dashed},
	]

\node { {$Q$}}
	child { node[below left] {$Q_1$}
		child { node {} }
		child { node[stencil point] {} }
		child { node {} }
		child { node[stencil point] {} }
	}
	child { node[below left] {$Q_2$}
		child { node[stencil point] {} }
		child { node[stencil point] {} }
		child { node[stencil point] {}
			child { node {$P$} }
			child { node {$P$} }
			child { node {$P$} }
			child { node {$P$} }
		}
		child { node[stencil point] {} }
	}
	child { node[below right] {$Q_3$}
		child { node {} }
		child { node[stencil point] {} }
		child { node {} }
		child { node {} }
	}
	child { node[below right] {$Q_4$}
		child { node[stencil point] {} }
		child { node[stencil point] {} }
		child { node {} }
		child { node {} }
	};

\end{scope}

\end{tikzpicture}}
\end{center}
\caption{Prediction stencil nodes to compute values at nodes $P$ within a given quadtree. The stencil comprises nodes belonging to four 
different quadtrees: $Q_i$, $i=1,2,3,4$.}\label{broth}
\end{figure}
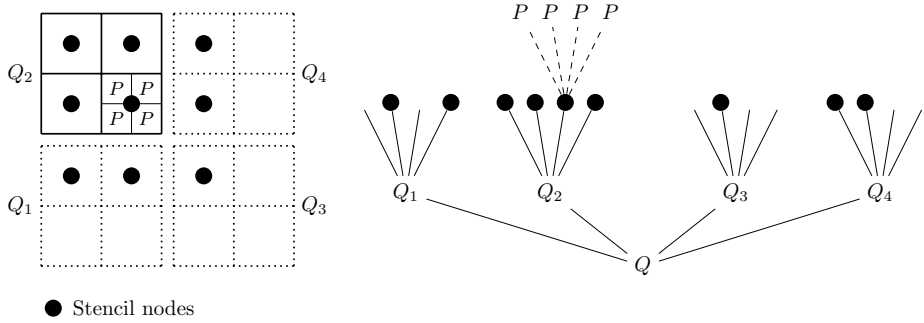

One of the main constraints on parallel implementations of the multiresolution algorithm 
has to do with the data locality required to efficiently perform the inter-level operations
throughout the tree data structure (quadtrees or octrees in two or three dimensions, respectively).
More specifically, the prediction operator needs to access all nodes in the interpolation stencil,
meaning that for any given node, one must access the parent-node at the successive coarser level 
and all its neighboring nodes, as shown in Figure \ref{broth} in the case of quadtrees. 
On the other hand, the projection operator needs to access the child-nodes of any given node
that is not a leaf.
Notice that it is very unlikely to achieve data locality in standard pointer-based tree data structures,
like the one considered in the \textsc{MBARETE} code.
In the latter case nodes are generated and linked recursively and even though pointers to the neighboring nodes are
introduced \cite{duarte:tel-00667857}, 
nodes belonging to a given interpolation stencil will hardly be stored contiguously in memory.
Additionally, other the projection and prediction operators, the refinement and coarsening operations, 
also performed on trees, should be executed in parallel as well, in order to maximize parallel coverage and minimize the impact of Amdahl's law.

The data structure  implemented in this work to support quadtrees and octrees, as well as the computation of the key multiresolution operations
are detailed in what follows.

\subsubsection{Data structure and implementation of trees}
\begin{figure}[h]
\begin{center}
\begin{tikzpicture}[scale=0.8]

\begin{scope}[xscale=10.0/48.0]
	\draw ( 0, 0.5) rectangle +(64, 1);
	\draw (48, 0.5) -- +(0, 1);
	\draw (52, 0.5) -- +(0, 1);

	\path
		( 0, 0) node (p0) {}
		(48, 0) node (p1) {}
		(52, 0) node (p2) {}
		(64, 0) node (p3) {};

	\begin{scope}[shorten >=0pt, shorten <=0pt, latex-latex]
		\scriptsize
		\draw (p0) edge node[below, align=center] {abscissa \\ (48 bits)} (p1);
		\draw (p1) edge node[below, align=center] {level \\ (4 bits)} (p2);
		\draw (p2) edge node[below, align=center] {free \\ (12 bits)} (p3);
	\end{scope}
\end{scope}

\end{tikzpicture}
\end{center}
\caption{Representation of a node using a $64$-bit integer.}\label{node}
\end{figure}
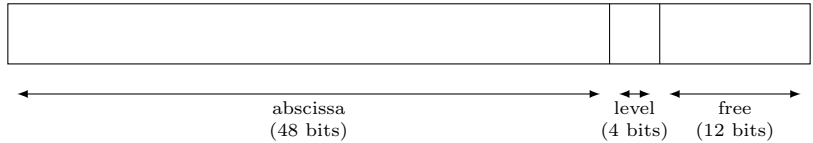

Without loss of generality we consider again the computational domain 
$\Omega=[0,1]^d$, $d=2,3$. 
The dyadic refinement of cells allows us to define the coordinates of a given node, $(x_1,\ldots,x_i)$, $i=1,d$,
using the binary numeral system.
That is, for each dimension $i$, the $x_i$-coordinate of a given node $P$ at the grid-level $j$ 
is given by $x_i=0.x_{i,1}x_{i,2}\cdots x_{i,j-1} x_{i,j}$,
where $x_{i,j-1}$ is equal to $0$ or $1$ depending on which of the two possible nodes
at grid-level $j-1$ generated the node $P$; and the same follows for $x_{i,j-2}$, up to $x_{i,1}$.
Using this convention one can easily define a space-filling curve \cite{MR1299533} throughout all nodes of the tree;
a powerful tool already employed in multiresolution
\cite{brix2009parallelisation,dahmen2011numerical} and AMR \cite{burstedde_p4est:_2011,weinzierl_peanotraversal_2011}
applications.
In our particular case 
we use a \emph{Morton order} space-filling curve \cite{morton66},
also known as \emph{Z-order} or \emph{Morton code}. For instance, in
two dimensions, for a given node $(x_1=0.x_{1,1}x_{1,2}\cdots
x_{1,j},\ x_2=0.x_{2,1}x_{2,2}\cdots x_{2,j})$, its Morton abscissa is constructed
by alternating the digits of $x_1$ and $x_2$:
$$s=x_{1,1}x_{2,1}x_{1,2}x_{2,2}\cdots x_{1,j}x_{2,j};$$ 
and the same follows in three dimensions.
In this way each node of the tree is uniquely defined by an integer.
For each node of the tree we can then store in a $64$-bit integer, its abscissa $s$
and its corresponding grid-level $j$.
In our implementation the first $48$ bits are intended to contain the abscissa and the following four, 
the grid-level, leaving the remaining bits available for tagging purposes during the multiresolution
operations, as shown in Figure \ref{node}.
In the case of three-dimensional problems, the latter choice allows one to encode nodes
of trees with up to $16$ grid-levels.

This representation of nodes 
is similar to the so-called CSAMR data structure introduced in \cite{ji_new_2010},
which then uses hash tables for efficient node lookup.
In our implementation we also employ hash tables, but over clustered nodes, in what follows denoted as \emph{blocks},
which enables further control on lookup granularity, and eases the operations for shared-memory parallelism.
Taking into account that $r= 0.s$ belongs to $[0,1[$, we partition $[0,1[$ into intervals.
The data structure is thus given by \emph{blocks} and \emph{collections of blocks}.
A block is then a structure that contains data related to a given interval, namely,
the interval bounds: $s_{\rm{min}}$ and  $s_{\rm{max}}$;
a vector containing all nodes
such that $s$ belongs to $[s_{\rm{min}}, s_{\rm{max}}[$;
and the size of the vector. 
A collection of blocks is given by pointers to
the blocks, ordered according to their $s_{\rm{min}}$ values. 
In practice, based on preliminary tests,
we operate with two collections of blocks, one 
for the leaf nodes and another one for the non-leaf nodes.

Abscissas are not necessary ordered inside 
blocks, and
binary search is used to locate the
block where the abscissa of a given node is possibly stored, followed by a
sequential search within that block. In our implementation the latter
strategy has proven to be the most efficient, in particular compared 
with one based on a fully ordered storage. A system of caches keeps
track of the last blocks accessed, effectively
decreasing the number of binary searches. These caches are local to TBB
tasks. 

Additionally, five auxiliary operations are implemented to properly handle this data structure.
Two of them aim at maintaining the size of the blocks within a predefined range,
by splitting into two a relatively large block or merging two neighboring blocks of relatively
small size.
The complexity of the search is hence of
$O(\log_2 N_b)$, where $N_b$ is the number of blocks in the
collection.
The next two operations are used to enlarge or reduce the size of the vectors
contained in the blocks.
These are performed by allocating a new vector of suitable size, copying the 
content of the original vector, and finally, deleting the original vector.
The fifth operation is implemented as a garbage collector, for instance, to effectively 
delete previously tagged nodes.
Notice that all these five operations are such that can be performed in parallel,
by using the \texttt{parallel\_for}
structure of the TBB library while looping over blocks.

Finally, the refinement process 
entails adding the abscissas of the new nodes at the end of the vector
of the block that contains the abscissa of the original node.
On the contrary, the coarsening process is performed using a \emph{lazy} approach;
that is, all nodes that need to be deleted are first tagged, using the free space
in their $64$-bit representation,
while in a second stage the garbage collector effectively remove them from their
corresponding blocks.
Further details are given in the following.

\subsubsection{Parallel implementation of the key multiresolution operations}\label{paraltree}

Grid adaptation in the multiresolution algorithm is essentially driven by the value
of details (\ref{eq:detail}), which in turn is computed using the projection and prediction operators.
In practice, 
details are computed in two stages, as follows.
\begin{enumerate}[I.]
 \item {\bf Projection.} Data is propagated from the leaf nodes to
  the root of the tree, an operation easily cast as a recursive computation. 
  Based on preliminary tests, we came out with the following implementation.
  Identifying a relatively low, fully populated level, $j_{\min}$, two recursive 
  computation are performed.
  First, for every node $P$ at grid-level $j_{\min}$, data is recursively transferred from 
  the corresponding leaves down to $P$. These are independent computations, hence performed
  in parallel by using recursive \texttt{parallel\_for} constructions.
  In a second step, a sequential computation propagates recursively data from level $j_{\min}$
  to the root of the tree.
  In this work, $j_{\min}$ was taken equal to $3$ and $5$, respectively, for three- and two-dimensional applications.
  It is worth noticing that often the sequential computational, although inexpensive, 
  needs not be performed, as in general no grid adaptation is present at those low grid-levels.

  \item {\bf Prediction and detail computation.} Details are computed for every node $P$ according to (\ref{eq:detail}),
  after applying the prediction operator by means of interpolation, in the following schematic manner.
\begin{algorithmic}
\For {$P \in \rm{\{Nodes\}}$}
 
   Define the list $L$ of Nodes in the prediction stencil associated with $P$ (see, \emph{e.g.}, Figure \ref{broth});
   
   Get values of Nodes in $L$;

   Compute the predicted value and the detail;
   
   Store the detail.

\EndFor
\end{algorithmic}
All these computations are easily implemented within a \texttt{parallel\_for} 
structure, looping over blocks.
\end{enumerate}

Once the details have been computed, the adapted grid is updated, also in parallel,
in two stages by means of the refinement and coarsening procedures, detailed in what follows.
Notice that in all cases a lazy strategy was again implemented.
\begin{enumerate}[I.]
 \item {\bf Refinement.} 
 Based on the detail values, nodes that need to be refined are first tagged in parallel.
 Then, all tagged cells are refined. The resulting tree will not necessary be a graded one.
 Therefore, additional nodes must be tagged for refinement in order to guarantee a 
 graded tree-structure. The process of tagging and refining is thus repeated until a graded
 tree-structure is achieved.
 
 The task-based refinement operation is implemented in the following way.
 Each available task considers a given set of consecutive blocks, 
 that is blocks of nodes within a given abscissa interval, say $[s_i,s_f[$.
 Whenever a node is created by refinement, it is directly stored in the corresponding 
 block that contains its abscissa, as previously explained.
 However, if necessary, two new blocks are created by each task to temporarily store those new nodes 
 whose abscissas are not contained in $[s_i,s_f[$ (one for abscissas smaller
 than $s_i$, and another one for those larger than or equal to $s_f$).
 These operations can be safely done in parallel because no task will be writing into any set or temporary block
 belonging to another task.
 Once all tagged nodes were refined, a new task-based operation starts, again 
 over consecutive sets of blocks, to move the previously created nodes from the temporary
 blocks to the corresponding block actually containing its abscissa.
 The latter process can be also safely performed in parallel because it only considers
 read-only accesses.
 
 \item {\bf Coarsening.} 
 At this point, leaf nodes that are no longer necessary can be eliminated.
 The latter applies to leaf nodes whose parent-nodes were not tagged for refinement
 in the previous stage,
 and that can be deleted without compromising the graded tree-structure.
 
 The task-based coarsening operation is implemented in the following way.
 Each available task considers a given set of consecutive blocks, and tags the leaf nodes
 that need to be deleted. Once the tagging phase is over, parent-nodes of tagged nodes are moved to the collection
 corresponding to leaf nodes, using the same task-based techniques described for the refinement.
 Finally, a parallel call to the garbage collector effectively deletes the tagged nodes.
 These three operations are repeated until no node deletion is needed. 
 Again these computations are safely performed in parallel using \texttt{parallel\_for} over blocks.
 
\end{enumerate}

\subsection{Parallel implementation of the splitting solver}

Before considering the parallel implementation itself, we have to define
the way of storing the unknowns of system (\ref{thesystem}) in memory.
This set of unknowns needs to be solved in every leaf node of the tree,
and there are basically two options to store them in memory,
following either an \emph{array of structures} or \emph{structure of arrays} approach.
The first option consists in storing all data related to a given node, in particular the unknowns, 
in a common structure, thus achieving node-wise memory contiguity.
On the other hand, the second option consists in storing each unknown in its own array,
thus achieving variable-wise memory contiguity.

In general,
within a splitting framework, the first layout turns out to be the most appropriate for the 
reaction solver, 
as it computes local reaction rates using all the unknowns at the current node.
The diffusion problem, on the other hand,  is solved for each unknown independently, 
hence a structure of arrays would provide better data locality.
Nevertheless, the reaction solver is highly compute-intensive and therefore, largely insensitive to the memory layout,
whereas the diffusion step displays very low compute-intensity.
Consequently, 
we have opted to store each unknown in its own contiguous array.
As we will later see, a large fraction of the runtime is spent in the reaction and diffusion steps;
having such a memory layout proves then 
to be very relevant, especially when the variables are not coupled in the diffusion operator.
The  Morton order 
defines the order of the unknowns, which are stored in vectors. 
Additionally,
the sparse matrices used during the diffusion steps can be stored according to a classical
\texttt{CSR} storage \cite{crs}.
However, an even more compact format was developed in this work
for diffusion problems with constant diffusion coefficients (see details in Appendix \ref{app:diffusion}).

As previously mentioned, 
the parallel implementation of the reaction solver is quite straightforward.
Once the access to the values of the unknowns $u$ at a given node is defined,
the reaction solver operates over the $N$ leaf nodes present in the adapted grid
as follows.
The reaction solver class is such that computes the solution for a given number of nodes, 
defined as a \texttt{blocked\_range<size\_t>} within the TBB framework.
Then, a call to \texttt{parallel\_for(blocked\_range<size\_t>(0,N),Reaction)}
parallelizes the computation, where \texttt{Reaction} corresponds to the solver class.

For the diffusion solver, two levels of parallelism is implemented.
Other the aforementioned parallel computation of the $m$ unknowns in system (\ref{thesystem}),
the matrix-vector products performed by \textsc{ROCK4} are also parallelized, 
using two \texttt{parallel\_for} constructions.
Further details on this \textsc{ROCK4} implementation is described in what follows

\subsubsection{Implementation of \textsc{ROCK4}}

\textsc{ROCK4} is an explicit time integration method, which can be viewed as the composition of
two methods that are successively performed \cite{MR1923724}. The first one is based on orthogonal polynomials and 
uses a three-term recurrence formula where the number of function evaluations, hence
stages of the method,
depends on the spectral radius of the Jacobian of the system.
The second one is an explicit four-stage Runge-Kutta scheme. For a system
$dy/dt=f(y)$ and a time step $\delta t$, the latter is given by
\begin{equation}\label{eq:RKscheme}
 \left.
 \begin{array}{ll}
  k_i= f\left(y_n+\delta t \ds \sum_{j=1}^{i-1} a_{i,j}k_j\right),& i=1,2,3,4,\\
  y_{n+1}= y_n+\delta t \ds \sum_{i=1}^4 b_i k_i,&
 \end{array}
\right\}
\end{equation}
where $a_{i,j}$ and $b_i$ correspond to the coefficients of the Runge-Kutta scheme.

Computing linear combinations of large vectors of \texttt{double} type,
as in (\ref{eq:RKscheme}),
is characterized by a rather low \emph{arithmetic intensity}\footnote{The 
arithmetic intensity is a measure of floating-point operations (FLOPs) 
performed by a given code (or code section) relative to the amount of memory 
accesses (Bytes) that are required to support those operations.},
less than $1/4$, because of the multiple memory accesses associated with 
reading and writing temporary arrays in memory. 
Nevertheless, 
achieving a high arithmetic intensity is fundamental to obtain
a high FLOP/s efficiency (see, \emph{e.g.}, \cite{arithI}).
In our implementation,
the arithmetic intensity is largely improved by exploiting the fact 
that an explicit Runge-Kutta scheme applied to a
linear problem $dy/dt=Ay$ is given by $y_{n+1}= \pi(\delta t A) y_n$, where
$\pi$ is a polynomial whose coefficients are a function of the $a_{i,j}$'s and
$k_i$'s coefficients \cite{MR1227985}, which are known and precomputed. 
In this way, we use the Horner's method to compute the four-stage
Runge-Kutta formula, thus reducing considerably the aforementioned memory
traffic.
Notice that a factorization of $\pi$ is not feasible
because this polynomial has complex roots. Finally,
the spectral radius $\rho(A)$ of $A$,
used to determine the number
of stages of the method,
is estimated 
using the Gershgorin circle theorem.

\section{Code performance and scalability}\label{sec:performance}
\begin{figure}[h]
\begin{center}
\includegraphics[width=0.7\linewidth]{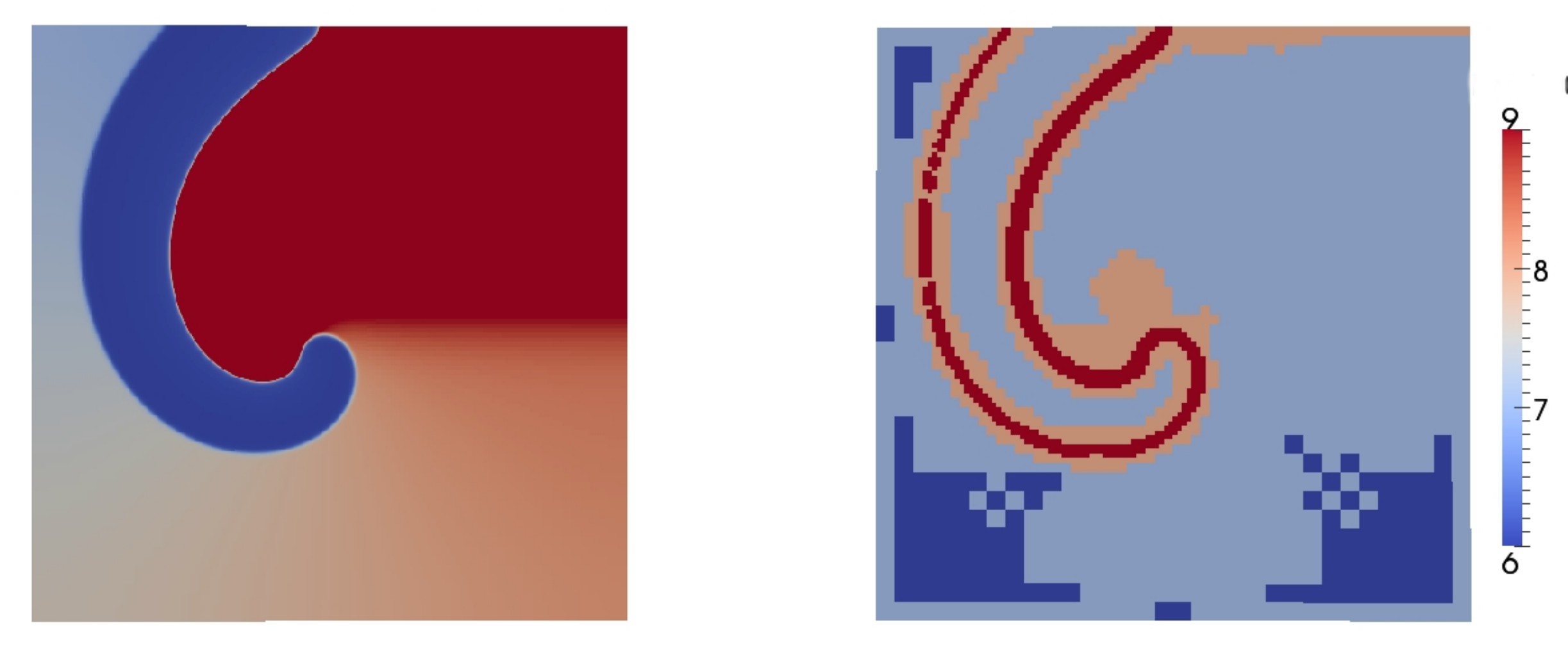}
\end{center}
\caption{BZ model, two-dimensional simulation over $9$ grid-levels.
Variable $u_1$ (left) with the corresponding adapted grid defined at different grid-levels (right).}
\label{bzwave2d}
\end{figure}

We have developed a C++ code, named
\texttt{Z-code}, which operates on multi-dimensional configurations;
the dimension $d$ is a template parameter for most classes and methods.
It runs successfully on modern Linux systems, and has been tested with both
\texttt{gcc} (version $4.8$ and higher) and the Intel C++ (version
$14.0.1$ and higher) compilers.
Figure \ref{bzwave2d} shows some results for the BZ problem.

Regarding the methods considered for the dedicated splitting solver,
codes for the \textsc{Radau5} and \textsc{ROCK4} methods are publicly
available on the Internet \cite{genevesoft}.
These codes, albeit very efficiently coded in Fortran 77, can
be substantially improved.
For \textsc{Radau5} we have carefully rewritten its implementation in C++,
thus removing a significant amount of
branching. The latter yields an enhanced
performance between $20$ to $30\,$\%. 
In the case of \textsc{ROCK4}, we have developed a completely new implementation,
efficiently adapted to linear problems,
that reduces memory accesses, as
described before\footnote{The sources of \texttt{Z-code} and related software can be provided
by contacting T.~Dumont, the corresponding author of the present
contribution. All codes are distributable under CeCILL-B licence 
(\url{http://www.cecill.info/licences.en.html}).}.

Numerical performance can be measured from many points of view.
In what follows we first present some
comparisons between computations on 
uniform Cartesian 
and multiresolution adapted meshes.
We then focus on the 
scalability attained for a given problem with increasing thread count (strong scaling).
The following tests were
performed on an Intel E5-2650 v3 platform (2 Haswell CPUs with 
$20$ cores in total and $40$ threads with hyperthreading)
and on a Xeon~Phi (MIC) 5110P ($60$ cores, $240$ threads
with simultaneous multithreading, and 8~GB of main memory)
using the Intel \texttt{icpc} compiler.

\subsection{Uniform Cartesian meshes versus multiresolution computations}

Grid adaptation results in significant data compression. 
However, it is worth verifying whether the latter effectively leads to performance enhancement
in terms of runtime,
given the additional computations required to dynamically adapt the grid
and the more elaborate implementation.

Nevertheless,
in order to guarantee a fair comparison, we have replaced the diffusion implementation in the 
\texttt{Z-code} with one more suited to uniform meshes, for all runs performed on a uniform grid.
Specifically, the matrix-vector products were replaced with a more efficient procedure
for five- and seven-point stencils, using two- and three-dimensional arrays, 
respectively, and cache blocking.
The latter avoids the additional cost of manipulating more complicated data structures and 
reduces the bandwidth of the computation.
The following results were obtained on the \textsl{Haswell} machine, using $40$
threads over $20$ cores. 

Recalling that $J$ stands for the maximum grid-level allowed for mesh refinement, 
Tables \ref{BZandKPP2D} and \ref{BZandKPP3D} gather the wall-clock computing times in seconds,
for various grid-levels $J$ and the aforementioned three problems in 
two and three dimensions.
Additionally, we can define a compression ratio in percentage, $CR$, as
$$
CR = \ds \left(1 - \frac{N}{N_J}\right)\times 100,
$$
where $N$ and $N_J$ correspond to the number of nodes in the adapted grid
and at the finest level, respectively.
In particular
$N_J$ is given by $2^{dJ}$ and is equal to the number of nodes in the corresponding uniform
Cartesian mesh.
During all these computations, the compression ratio remained between $80\,$\% and $83\,$\%.
\begin{table}[h]
\begin{center}
\begin{tabular}{|l|l||c|c|c|}
\cline{3-5}
\multicolumn{2}{c}{}
\tlvs
&\multicolumn{3}{|c|}{$J$}\\
\cline{3-5}
\multicolumn{2}{c|}{}
\tlvs
&8&9&10\\
\hline
\tlvs
\multirow{2}{*}{\textbf{NAGUMO}}
&MR&$\boldsymbol{2.51\times 10^{-3}}$ (1.7)&  $3.02\times 10^{-3}$ \phantom{(1.0)}&$4.97\times 10^{-3}$ \phantom{(1.0)}\\
&CM&  $1.49\times 10^{-3}$ \phantom{(1.0)}&  $\boldsymbol{8.27\times 10^{-3}}$ (2.7)&$\boldsymbol{2.09\times 10^{-2}}$ (4.2)\\
\hline
\hline
\tlvs
\multirow{2}{*}{\textbf{BZ}}&MR&$1.14\times 10^{-2}$ \phantom{(1.0)}&  $2.51\times 10^{-2}$ \phantom{(1.0)}&$4.06\times 10^{-2}$ \phantom{(1.0)}\\
&CM&  $\boldsymbol{1.31\times 10^{-2}}$ (1.2)&  $\boldsymbol{4.06\times 10^{-2}}$ (1.6)&$\boldsymbol{1.78\times 10^{-1}}$ (4.4)\\
\hline
\hline
\tlvs
\multirow{2}{*}{\textbf{STROKE}}&MR&$2.4\times 10^{-2}$ \phantom{(10.0)} &  $4.03 \times 10^{-2}$ \phantom{10.0} & $1.03 \times 10^{-1} $ \phantom{(10.0)}\\
&CM&  $\boldsymbol{3.2 \times 10^{-1}}$ (13.3)&  $\boldsymbol{1.21}$\phantom{$\times 10^{-1}$}  (30.0)&  $\boldsymbol{4.80}$\phantom{$\times 10^{-1}$} (46.6)\\
\hline
\end{tabular}
\end{center}
\caption{Wall-clock times in seconds for two-dimensional multiresolution (MR) and 
uniform Cartesian mesh (CM) computations. Largest values are indicated in bold with the
corresponding ratio in parenthesis.}\label{BZandKPP2D}
\end{table}
\begin{table}[h]
\begin{center}
\begin{tabular}{|l|l||c|c|}
\cline{3-4}
\multicolumn{2}{c}{}
\tlvs
&\multicolumn{2}{|c|}{$J$}\\
\cline{3-4}
\multicolumn{2}{c|}{}
\tlvs
&8&9\\
\hline
\tlvs
\multirow{2}{*}{\textbf{NAGUMO}}
&MR& $0.13$ \phantom{(1.0)}&$0.31$ \phantom{(1.0)}\\
&CM& $\boldsymbol{0.54}$ (2.5)&$\boldsymbol{2.47}$ (4.5)\\
\hline
\hline
\tlvs
\multirow{2}{*}{\textbf{BZ}}&MR&$0.97$ \phantom{(1.0)}&\phantom{0}$5.75$ \phantom{(1.0)}\\
&CM& $\boldsymbol{3.05}$ (3.1)&$\boldsymbol{23.60}$ (4.1)\\
\hline
\hline
\tlvs
\multirow{2}{*}{\textbf{STROKE}}&MR&\phantom{0} $1.53$ \phantom{(10.0)}&\phantom{0} $10.81$ \phantom{(10.0)}\\
&CM& $\boldsymbol{76.68}$ (50.0)&$\boldsymbol{610.50}$ (56.5)\\
\hline
\end{tabular}
\end{center}
\caption{Wall-clock times in seconds for three-dimensional multiresolution (MR) and 
uniform Cartesian mesh (CM) computations. Largest values are indicated in bold with the
corresponding ratio in parenthesis.}\label{BZandKPP3D}
\end{table}

These results clearly highlight the advantage of using the multiresolution
algorithm to solve stiff problems disclosing localized reaction fronts
with an appropriate grid resolution.
As a matter of fact, looking at the previous results we see that
computations on a uniform grid are more efficient only for the 
two-dimensional, non-stiff NAGUMO problem on a rather coarse grid.
In general 
the additional cost related to mesh adaptation and more complex data structures
is counterbalanced by the computing savings achieved with a compressed representation.
The latter is even more relevant for problems with very stiff reaction terms,
confirming the findings and forecasts previously described in Appendix \ref{fixedmeshes}.

\subsection{Analysis of scalability}
We investigate the strong scaling behavior of \texttt{Z-code} 
by timing one single simulation time step, and breaking down the total step time (S) into
the mesh adaptation (A), reaction solver (R), and diffusion solver (D).
Notice that the adapted grid is updated at every time step,
the detail computation is hence included in the mesh adaptation process,
and the diffusion solver must also construct the matrices used in the linear systems.

The three reference problems probe different regimes of compute-intensity.
The NAGUMO problem is characterized by a rather inexpensive, non-stiff reaction term,
whereas the computing effort for
the BZ problem is roughly balanced between the reaction solver and the other two procedures;
that for the STROKE problem, on the other hand, is completely dominated by a very expensive, stiff reaction term.
Table \ref{time_percent}
shows the percentages of wall-clock time spent on each of the three major code phases,
for three-dimensional simulations performed over $10$ grid-levels ($J=10$) and $40$ threads on the Haswell platform.
\begin{table}[h]
\begin{center}
\begin{tabular}{|l|c|c|c|}
\cline{2-4}
\multicolumn{1}{c|}{}
\tlvs
& \textbf{NAGUMO}
& \textbf{BZ}& \textbf{STROKE}\\
\hline
\tlvs
Mesh Adaptation 
&\textbf{68.00}&19.23&\phantom{0}7.20\\
\hline
\tlvs
Reaction Solver 
&\phantom{0}0.92&\textbf{56.23}&\textbf{88.89}\\
\hline
\tlvs
Diffusion Solver 
&31.08&24.54&\phantom{0}3.91\\
\hline
\end{tabular}
\end{center}
\caption{Percentages of wall-clock time for three-dimensional simulations
with maximum grid-level $J=10$ and $40$ threads on the Haswell platform.
Dominant contributions are indicated in bold.}\label{time_percent}
\end{table}

\begin{figure}[!h]
\centering
\begin{tabular}{cc}
\includegraphics[width=0.3385\textwidth]{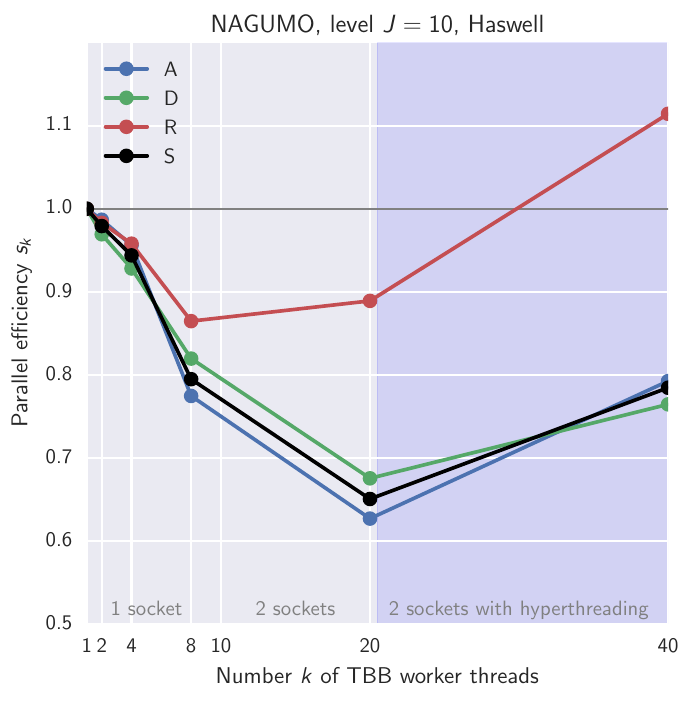} &
\includegraphics[width=0.3385\textwidth]{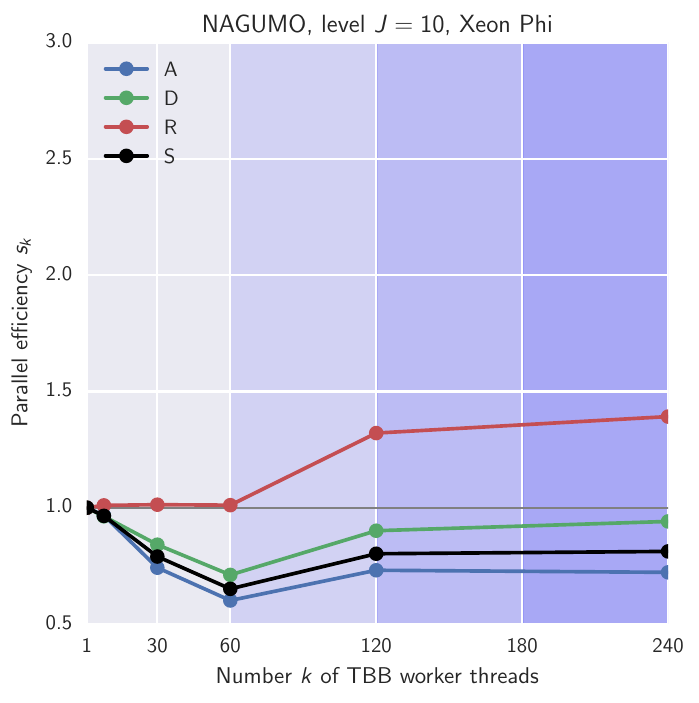} \\
\includegraphics[width=0.3385\textwidth]{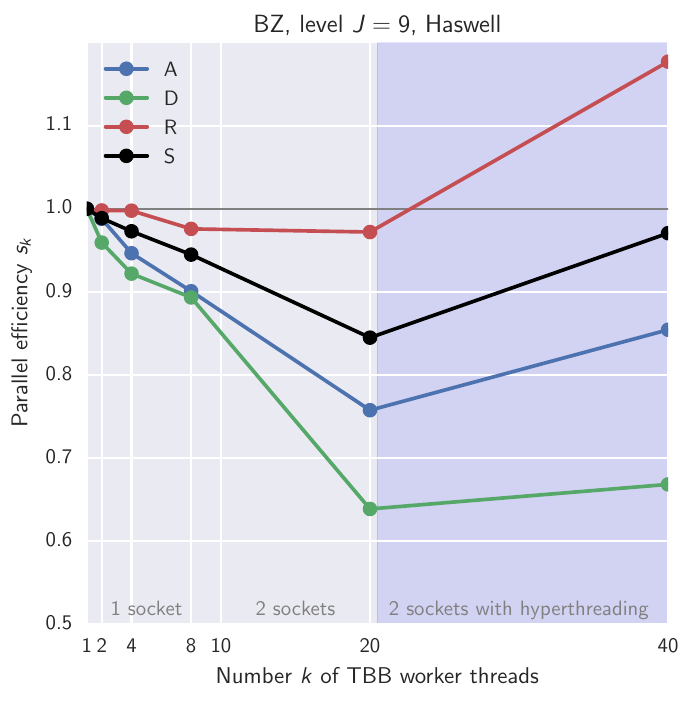} &
\includegraphics[width=0.3385\textwidth]{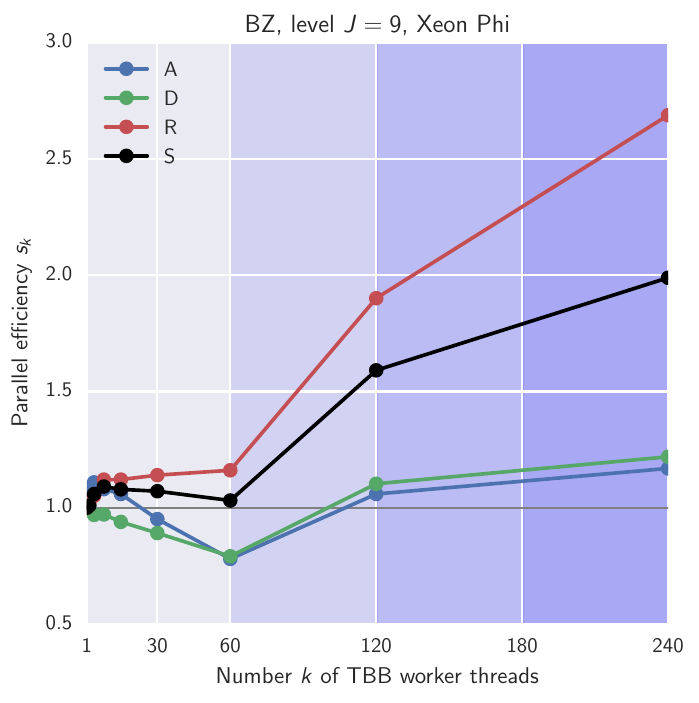} \\
\includegraphics[width=0.3385\textwidth]{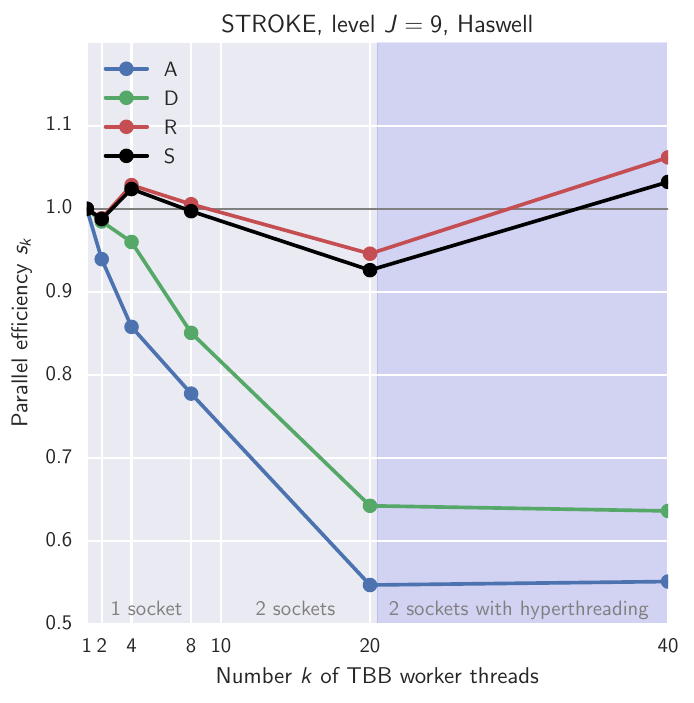} &
\includegraphics[width=0.3385\textwidth]{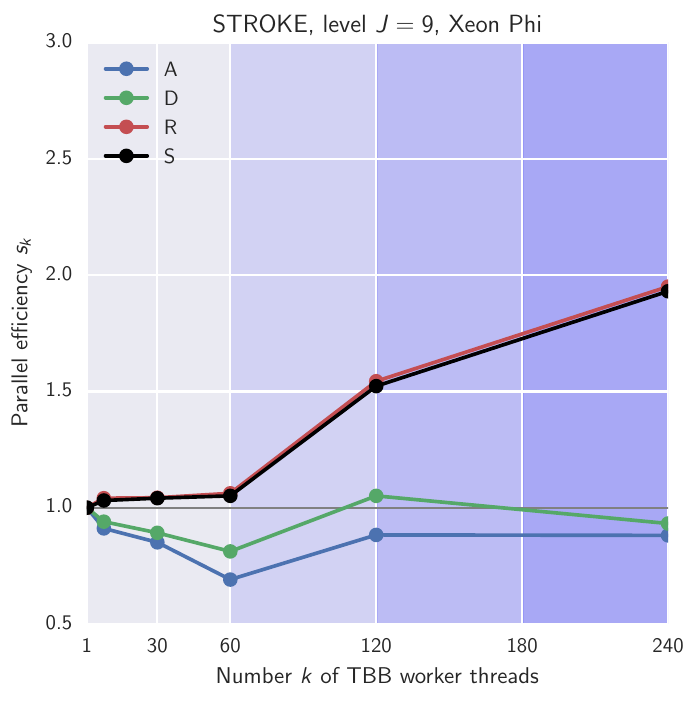}
\end{tabular}
\caption{Parallel efficiency $s_k$ in terms of number $k$ of TBB worker threads for the NAGUMO (top), BZ (middle), and STROKE (bottom) problems on Haswell (left) and Xeon~Phi (right) architectures, for mesh adaptation (A), diffusion (D), reaction (R), and total time step (S).
The blue-shaded areas highlight thread counts with more than one TBB worker thread per CPU core (2 threads/core on Xeon Haswell and 2--4 threads/core on Xeon Phi).}\label{s_k}
\end{figure}

We now discuss scalability in terms of parallel efficiency.
Let $T_k$ be the wall-clock computing time in seconds when running on $k$ threads. We define the
parallel efficiency on $k$ threads by
$$
s_k=\frac{T_1}{n_k T_k},
$$ 
where $n_k$ corresponds to the number of cores used when computing over $k$ threads.
In particular 
because of the \emph{simultaneous multithreading}, 
that is, running more than one thread per core,
$n_k$ is taken either as
$$n_k=\min(k,20),$$ 
for the Haswell platform, or
$$n_k=\min(k,60),$$
for the Xeon~Phi.
In this way  
we compute the efficiency based on the number of CPU \emph{cores},
not the number of available hardware threads.
The motivation for this choice is that cores, unlike hardware threads,
constitute independent compute units.
The scalability plots therefore present scaling results as a function of independent hardware resources, \emph{i.e.} cores, participating in the computation.
An important consequence is that, when using hyperthreading, efficiencies may well exceed one.
This simply means that running more than one thread per core may result in a more efficient execution at the core level
than what can be achieved with a single thread.
We will discuss hyperthreading in more detail in the following.
Figure~\ref{s_k} illustrates the parallel efficiency, $s_k$,
for all three problems, on Haswell and Xeon Phi, in three-dimensional simulations.

Overall, the code achieves good shared-memory scalability on compute-intensive problems (BZ and STROKE), with efficiencies over $80\,$\% over 2 Haswell sockets.
The scaling of the reaction step is very good mainly because of its embarrassingly parallel nature and the dynamic load balancing provided by TBB; there is actually no shortage of parallelism in this part of the time step.
In addition, since the reaction step is strongly compute bound, cores do not compete for shared resources, allowing for linear scaling.
Only for the NAGUMO problem, whose reaction computation is quite inexpensive, does the reaction efficiency drop below $90\,$\%.

The diffusion and adaptation steps achieve lower parallel efficiencies than the reaction one.
A first reason for this is that both phases are more memory-intensive than the reaction; therefore, more pressure is put on the memory controllers, which are shared between cores of a given socket.
This results in less-than-ideal scaling, as memory accesses become a bottleneck.
Note that for the diffusion operator, we have used the
compact structure described in Appendix \ref{app:diffusion},
which effectively improved performance over a classical \texttt{CSR} structure for the problems studied here.
Other the considerations on memory access, parallelization of the diffusion solver and multiresolution operations is intrinsically more complex than that of the reaction solver, yielding slightly less exposed parallelism and a degradation of parallel performance due to Amdahl's law.
This highlights the importance of exposing as much parallelism as possible, in particular in complicated algorithms such as mesh adaptation via multiresolution.

We now discuss the impact of simultaneous multithreading (hyperthreading) on the performance of \texttt{Z-Code} on the Haswell Xeon and Xeon Phi platforms.
Hyperthreading allows the cores' execution units to be used by a second application thread whenever the core stalls executing an instruction stream from the first thread.
Typically, it is beneficial for applications that have good thread scalability, and where the execution is frequently stalled by latency-related events such as cache misses, branches, or bad pipelining.

For \texttt{Z-Code} on Haswell, moving from 1 to 2 threads per core (20 to 40 TBB worker threads) provides significant performance improvements, especially in the reaction solver where hyperthreading gains up to 1.25 speedup, with a geomean across problems of 1.20.
The reaction solver features a very complex control flow in the \textsc{Radau5} algorithm,
together with BLAS function calls and operations on small matrices.
All these contribute to core stalls due to bad instruction pipelining and complex branching.
These latencies, together with the good thread scalability of the reaction solver, explain the significant performance gains attained by hyperthreading.

By contrast, still on Haswell, the gains are more modest for the multiresolution algorithm (geomean 1.13) and even more so for the diffusion solver (geomean 1.05).
The diffusion step features simpler control flow and, as previously discussed,
its performance is dominated by memory controller bottlenecks,
which hyperthreading cannot alleviate because they lie outside of the CPU core.

In regard to the Xeon Phi, on our Knights Corner coprocessor, simultaneous multithreading is crucial to enhancing performance.
Specifically, it allows hiding some of the execution latencies that cannot be overlapped by the in-order architecture.
In addition, because the Xeon Phi features two pipelines per core, it cannot reach its peak performance at less than 2 threads per core.
Just like on the Haswell architecture, we find the largest performance gains on the reaction solver, with a geomean speedup of 1.80 across problems between 60 and 240 threads (1 and 4 threads/core).
On Xeon Phi, however, hyperthreading is also beneficial to multiresolution mesh adaptation (geomean 1.32) and diffusion (geomean 1.33), highlighting the importance of also overlapping the latencies of memory accesses.

\section{Concluding remarks and prospects} 

The main purpose of the present work was to introduce shared-memory parallelism in a thoughtful way,
to a tailored numerical strategy previously developed to tackle the simulation of stiff reaction-diffusion models
disclosing propagating fronts \cite{MR2890259}.
In order to achieve this,
a new data structure was conceived along with parallel-friendly implementations of the key algorithms into play,
including the \textsc{ROCK4} and \textsc{Radau5} solvers.
Shared-memory parallelism was introduced using the task-based TBB runtime library, 
which turned out to be particularly well suited to our purposes.

The numerical performance of the resulting implementation was assessed for three reaction-diffusion models,
reaching a very satisfactory level of efficiency on shared-memory architectures such as the Intel E5-2650 
and the Xeon~Phi (MIC) 5110P.
For models involving a certain level of modeling complexity such as BZ or STROKE, 
high scalability is achieved mainly because of the very efficient performance of the
reaction solver within a splitting context.
In addition to that, the present parallel implementation of the diffusion solver and the multiresolution algorithm 
shows also to perform reasonably well, given the limitations commonly encountered in standard computing architectures,
as previously discussed.
Finally, simultaneous multithreading is shown to improve considerably the computational performance, in particular on many-core architectures.

The proposed strategy is thus suitable for the simulation of multi-dimensional time-space multi-scale problems, 
which would be out of reach using classical approaches on standard computing resources.
In particular 
these reaction-diffusion problems are representative building blocks of more complex and realistic applications,
encountered, for instance, in biomedical engineering \cite{Dumont2013}, combustion \cite{duarte:hal-00727442}, or plasma physics \cite{MR2869398}. 
In this context ongoing work encompasses the extension to arbitrary domains, taking into account, for example, 
that the STROKE model applies to the human brain \cite{Dumont2013};
and to more general diffusion problems, as used, for example, in multi-component transport modeling for combustion problems \cite{MR3213800}.
A key extension of the present work deals with the implementation of hybrid parallelism using
both shared- and distributed-memory formalisms.
In this regard we believe the implementations and findings of this work are relevant to distributed-memory codes as well.
Achieving strong scaling is particularly challenging for multiresolution and AMR codes,
because of load imbalance and fine-grain communications.
For some applications, shared-memory parallelism with dynamic load balancing
may allow distributing work evenly between the CPU cores within an MPI rank,
thereby helping reduce total imbalance.
We thus expect shared-memory parallelism to gain in relevance for complex and imbalanced MPI applications.
These are topics of our current research.

\section*{Acknowledgments}
This development was funded by  grants from  the ANR project (French National Research Agency - ANR Blancs)
S\' echelles (2009-2013) and from the French \emph{AMIES} \emph{Peps} program.
Part of this work was conducted by the Exascale Computing Research laboratory,
thanks to the support of CEA, GENCI, Intel, UVSQ.
Any opinions, findings, and conclusions or recommendations expressed in this material are those of the author(s) and do not necessarily reflect the views of CEA, GENCI, Intel or UVSQ.

\appendix

\section{Numerical experimentation on uniform Cartesian meshes}\label{fixedmeshes}  
Let us consider the numerical performance of the dedicated Strang scheme 
introduced in \S\ref{sec:splitting},
to solve (\ref{thesystem}) on a uniform Cartesian grid.
Table \ref{diff_percent} shows the percentages of total wall-clock time spent in the diffusion step 
of the splitting strategy.
\begin{table}[h]
\begin{center}
\begin{tabular}{|c|c|c|c|}
\hline
\tlvs
\textbf{Grid size} & \textbf{NAGUMO}& \textbf{BZ}& \textbf{STROKE}\\
\hline
\tlvs
$512^3$&$92.5$ &$34.5$&$1.6$\\
\hline
\tlvs
$1024^2$&$ 88.5$ &$20.1$&$1.0$ \\
\hline
\end{tabular}
\end{center}
\caption{Percentages of wall-clock time employed to solve the diffusion problem.}\label{diff_percent}
\end{table}

Notice that for the models of higher complexity, that is,
BZ and STROKE, most of the simulation time is spent in the reaction step.
Consequently, as a measure of the computational complexity, we use
the number of evaluations of the right-hand side performed by the
\textsc{Radau5} program on a given grid-node during the $R_{\Delta t/2}$-step. 
This complexity measure may considerably vary throughout the computational  
domain, given that the splitting time step $\Delta t/2$
can be further split according to the time-stepping procedure
and that an iterative, simplified Newton method is implemented
to solve the nonlinear systems.
In particular the minimum complexity corresponds to grid-nodes where
only one reaction sub-step is required during a given $\Delta t/2$, 
as well as one iteration of the Newton method.
\begin{figure}[h]
\begin{center}
\raisebox{1ex}{\includegraphics[width=0.4\linewidth]{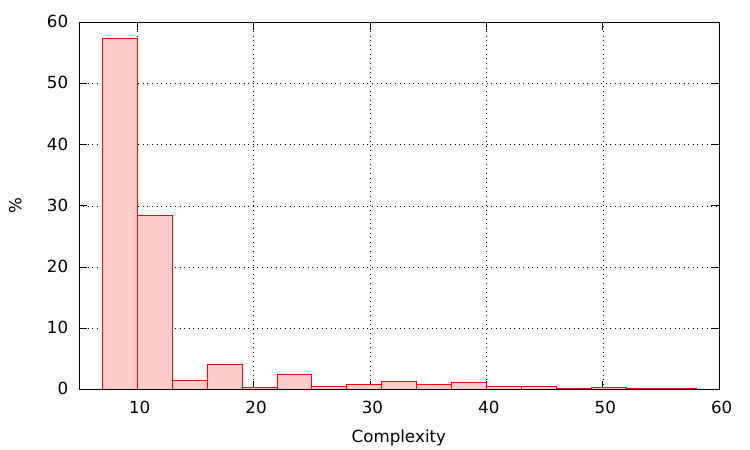}}
\includegraphics[width=0.25\linewidth]{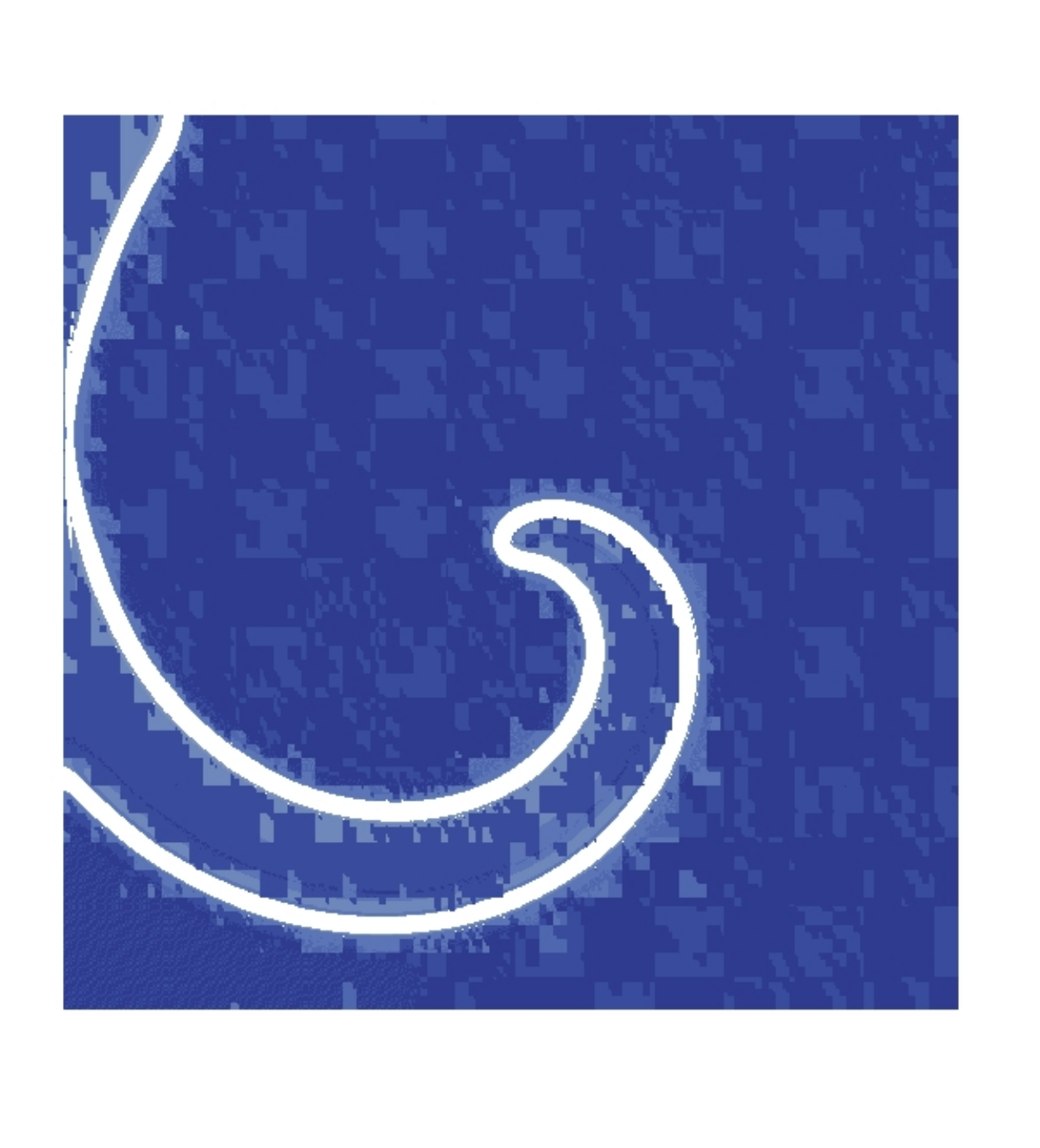}
\end{center}
\caption{BZ problem. Histogram of computational complexity (left). In the right
  picture, white zones correspond to points where the complexity is
  greater than $17$.}
\label{cplxbz}
\end{figure}
\begin{figure}[h]
\begin{center}
\includegraphics[width=0.4\linewidth]{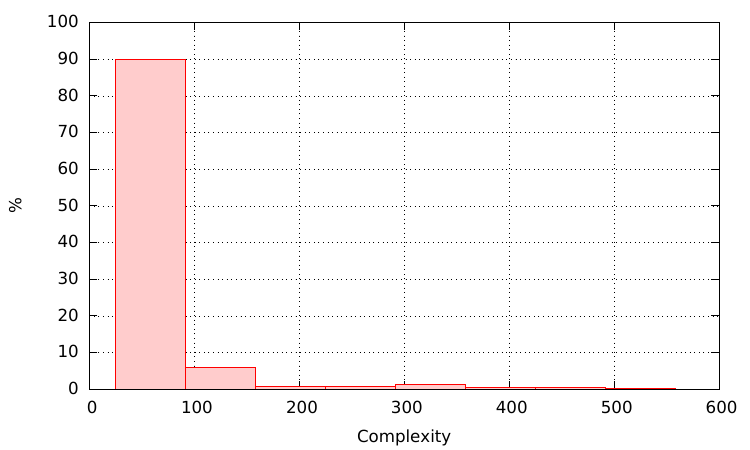}
\raisebox{2ex}{\includegraphics[width=0.22\linewidth]{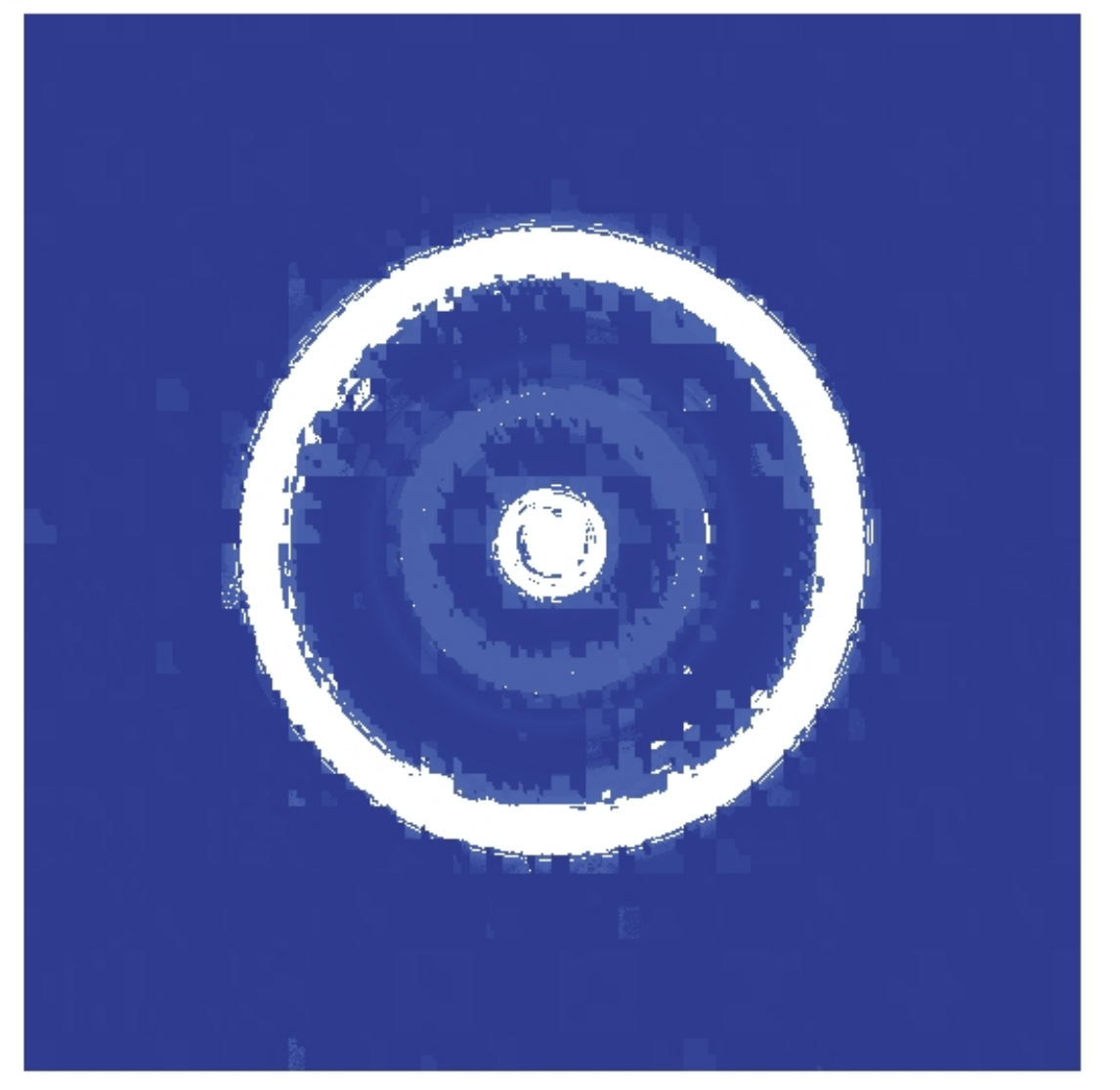}}
\end{center}
\caption{STROKE problem. Histogram of computational complexity (left). In the right
  picture, white zones correspond to points where the complexity is
  greater than $80$.}
\label{cplxstroke}
\end{figure}

Considering a uniform mesh of $1024^2$ grid-nodes, 
Figures \ref{cplxbz} and \ref{cplxstroke} show the computational complexity
measured for a given $\Delta t/2$ for the BZ and STROKE models, respectively.
For the BZ problem
the Jacobian is computed analytically involving three right-hand side evaluations.
In this case
$87\,$\% of the CPU time is spent over regions where the complexity is less than $17$. 
These are the regions where the solution is close to
the reaction equilibrium and therefore, the complexity remains also close to its
minimum value.
Similarly, about $90\,$\% of the reaction CPU time 
is spent in regions where the complexity is less than $60$ 
for the STROKE model.
Here the Jacobian is computed numerically.

\section{Choice of TBB for multiresolution applications}\label{TBBvsOpenMP}
In this work we have initially focused our attention on OpenMP API and TBB, both well-established options for shared-memory parallelism.

The OpenMP standard is a widely supported and popular choice, especially for high-performance computing (HPC) codes.
It relies on compiler directives to introduce parallel regions and
work-sharing constructs and tasks, which are executed in parallel by a runtime library component.
OpenMP has been supporting tasking since version 3, and has recently introduced more sophisticated task constructs with version 4,
with the \texttt{taskgroup} and \texttt{depend} constructs.

A key design characteristic of the OpenMP API is the reliance on a thread-based approach.
That is, parallelism is defined explicitly by parallel regions which enforce parallel execution,
and work and tasks are scheduled across threads participating in their containing parallel region.
The set of threads taking part in a region is determined when entering the region, and fixed for its entire duration.
Unfortunately, this severely limits opportunities for exploiting nested parallelism with OpenMP, even though the standard supports nested parallel regions.
Actually, due to the fact that threads are assigned to a region for its whole duration, there can be no dynamic load balancing or work-sharing \emph{across} parallel regions.
In addition, creating an OpenMP region requires knowing about the outside parallel execution context.
That is, since OpenMP parallelism is mandatory, nesting parallel regions may lead to exponential creation of OpenMP worker threads for some runtime implementations,
resulting in oversubscription.
Finally, since all threads of a parallel region are required to enter an OpenMP work-sharing construct,
one cannot arbitrarily nest OpenMP tasks and work-sharing constructs.
In particular one cannot introduce parallelism by nesting \texttt{omp for} constructs within \texttt{task} constructs.

Because of the aforementioned points, OpenMP parallelism is not \emph{composable}, that is, 
it is not possible to expose parallelism in an opportunistic way,
regardless of the outer context from which a piece of code may be called.
For the development of a C++ multiresolution code, we found composability to be a desirable property,
as it allows introducing parallelism at any level in the code, without knowledge or constraint on the outside calling context.
This is particularly important for codes with complex call graphs and multiple call paths, as often found in C++ when using classes for code reuse.

Intel TBB, on the other hand, relies on a different approach to parallel programming.
It is an open-source C++ template library for task parallelism, and requires no special support in the compiler.
In addition to task-based parallel constructs, TBB also provides concurrent data structures and memory allocators.
The library implements a composable task-based runtime, meaning that
any C++ function can expose parallelism with (potentially arbitrarily nested) constructs such as \texttt{parallel\_for} or \texttt{parallel\_reduce},
irrespective of the calling context.
Dynamic load balancing is achieved through a work-stealing runtime.
Since both tasks and constructs such as \texttt{parallel\_for} rely on the same task-based scheduler, 
TBB parallel idioms may be combined and nested arbitrarily,
as long as the programmer avoids races through appropriate synchronizations.
In particular, TBB's composability allows us to achieve good scalability with the adaptive multiresolution algorithm (see \S~\ref{paraltree}) by recursively nesting \texttt{parallel\_for} constructs to expose parallelism in a way that follows the tree structure.
In addition to a composable task runtime, TBB provides a parallel algorithm library
with support for parallel reduction and sorting.

We chose TBB over an OpenMP implementation mainly because the former provides
an efficient task runtime with dynamic load balancing, composability of parallel codes,
and a set of parallel data structures and algorithms.
Equivalent functionality could have been achieved using OpenMP tasks,
but more effort would have been necessary to circumvent some of the aforementioned shortcomings.

\section{A compact data structure for diffusion sparse matrices}\label{app:diffusion}
When the diffusion coefficients, $\varepsilon_i(x)$ in
(\ref{thesystem}), are constant, it is possible to further reduce the
amount of memory necessary to store the matrix entries associated with the
diffusion problem, leading to lower memory bandwidth requirements for the program.
Notice that for a given node, the matrix entries defined by the spatial 
discretization of the Laplace operator can have three possible values, 
depending on whether the neighboring
nodes belong to the same, upper or lower grid-level.
We therefore use the data structure shown in Figure \ref{newdata}.
It consists of two arrays;
the array \textit{A} contains pointers to the beginning of lines in the second array $K$,
where each line corresponds to one node in the adapted grid.
The array \textit{K} contains \emph{line descriptors}, 
one for each node in the adapted grid, defined as follows. 
The first integer of a line descriptor is a 32-bit integer, containing the grid-level $j$ of the current 
node, followed by the number of neighboring nodes at grid-level $j$, $j+1$, and $j-1$,
as shown in Figure \ref{newdata}. The following integers then indicate the corresponding
neighboring nodes.
As an illustration, Figure \ref{newdata1} shows a typical line descriptor.
\begin{figure}[h]
\centering
\scalebox{0.8}{
\begin{tikzpicture}
\usetikzlibrary{calc}
\begin{scope}[scale=0.5, yshift=2cm]
	\draw (0, 0) rectangle +(10, 1);
	\draw (2, 0) -- +(0, 1);
	\draw (3, 0) -- +(0, 1);
	\draw (4, 0) -- +(0, 1);
	\draw (0, 0) node[above left] {$A$};
	\draw (2.5, 1) node[above] {\tiny $i$};
	\draw (3.5, 1) node[above] {\tiny $i$+1};
	\coordinate[anchor=south] (i) at (2.5, 0);
	\coordinate[anchor=south] (f) at (3.5, 0);
\end{scope}

\begin{scope}[scale=0.5, yshift=0cm]
	\draw (0, 0) rectangle +(13, 1);
	\draw (3.5, 0) -- +(0, 1);
	\draw (6.5, 0) -- +(0, 1);
	\node[above left] (0, 0) {$K$};
	\node[above] at (5, 0) {\scriptsize Line $i$};
	\coordinate[anchor=north] (i2) at (3.6, 1);
	\coordinate[anchor=north] (f2) at (6.6, 1);
\end{scope}

\draw (i) edge[gray, ->, out=-90, in=90] (i2);
\draw (f) edge[gray, ->, out=-90, in=90] (f2);

\begin{scope}[xshift=7cm, xscale=0.4, yscale=0.5]
\small
\draw ( 0, 0) rectangle (20, 1) ++(0, -1) coordinate (p3);
\draw ( 1, 0) -- +(0, 1);
\draw ( 2, 0) -- +(0, 1);
\draw ( 3, 0) -- +(0, 1);
\draw ( 4, 0) coordinate (p0) -- +(0, 1);
\draw (10, 0) coordinate (p1) -- +(0, 1);
\draw (16, 0) coordinate (p2) -- +(0, 1);
\small
\node at (0.5, 0.5) {$j$};
\node at (1.5, 0.5) {$n_1$};
\node at (2.5, 0.5) {$n_2$};
\node at (3.5, 0.5) {$n_3$};
\end{scope}

\begin{scope}
\draw ($(p0)+(0,-4mm)$) edge[latex-latex] node[below, align=center] {$n_1$ ints} ($(p1)+(0,-4mm)$);
\draw ($(p1)+(0,-4mm)$) edge[latex-latex] node[below, align=center] {$n_2$ ints} ($(p2)+(0,-4mm)$);
\draw ($(p2)+(0,-4mm)$) edge[latex-latex] node[below, align=center] {$n_3$ ints} ($(p3)+(0,-4mm)$);
\end{scope}

\end{tikzpicture}
}
\caption{Compact data structure for diffusion matrices. The right part shows the structure of a line descriptor.} 
\label{newdata}
\end{figure}
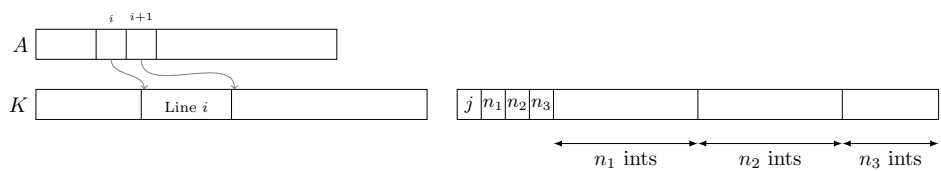
\begin{figure}[ht]
\centering
\subfloat
{
\includegraphics[width=.2\linewidth]{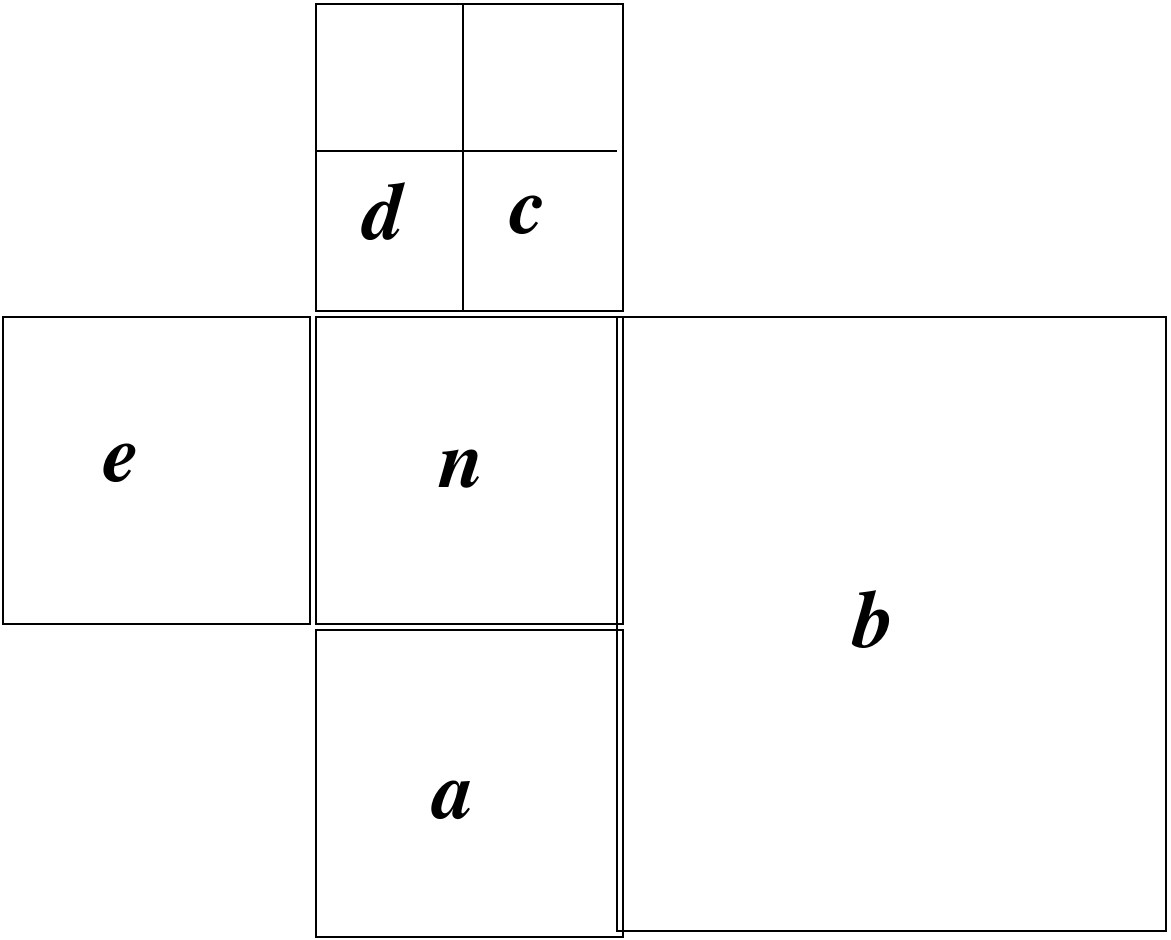}
}
\quad
\subfloat
{
\includegraphics[width=.4\linewidth]{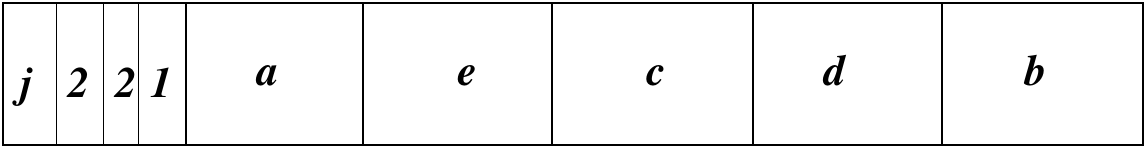}
}
\caption{Illustration of a line descriptor for a given node $n$ at grid-level $j$.}
\label{newdata1}
\end{figure}
\begin{table}[ht]
\begin{center}
\begin{tabular}{|l|c|c|c|c| }
\cline{2-5}
\multicolumn{1}{c}{}
\tlvs
& \multicolumn{2}{|c|}{Haswell} &
\multicolumn{2}{|c|}{Xeon Phi}\\
\cline{2-5}
\multicolumn{1}{c|}{}
\tlvs
& $d=2$, $J=10$& $d=3$, $J=10$&$d=2$, $J=10$& $d=3$, $J=9$\\ 
\hline 
\tlvs
\textbf{Ratio} &0.98&0.67&0.97&0.86\\
\hline
\end{tabular}
\caption{Ratio between wall-clock times using the compact and \texttt{CSR} storages.
BZ model simulated in dimension $d$ and with maximum grid-level $J$}
\label{cpmat}
\end{center}
\end{table}

Table \ref{cpmat} compares the wall-clock times obtained when computing
one diffusion step with the classical \texttt{CSR} data structure and 
the aforementioned compact storage, for the BZ model.
Simulations on the Haswell platform were performed using $20$ cores and $40$ threads, 
whereas for Xeon Phi, $60$ cores and $240$ threads were considered. 
Computational gains are rather weak for two-dimensional simulations, for 
which the matrices remain relatively small.
However, better performances are obtained with the compact structure
when dealing with three-dimensional matrices.

\bibliography{u}
\bibliographystyle{plain}

\end{document}